\documentclass{birkjour}

\usepackage{stmaryrd}
\usepackage{amsfonts}
\usepackage{bbm}
\usepackage{amscd}
\usepackage{mathrsfs}
\usepackage{latexsym,amssymb,amsmath,amscd,amscd,amsthm,amsxtra,xypic}

\newtheorem{thm}{Theorem}[section]
\newtheorem{lem}[thm]{Lemma}
\newtheorem{cor}[thm]{Corollary}
\newtheorem{pro}[thm]{Proposition}
\newtheorem{ex}[thm]{Example}
\newtheorem{rmk}[thm]{Remark}

\newcommand{\lon }{\,\rightarrow\,}
\newcommand{\be }{\begin{eqnarray*}}
\newcommand{\ee }{\end{eqnarray*}}

\newcommand{\defbe}{\triangleq}
\newcommand{\pf}{\noindent{\bf Proof.}\ }

\newcommand{\Real}{\mathbb R}

\newcommand{\huaF}{\mathcal{F}}
\newcommand{\huaG}{\mathcal{G}}

\newcommand{\CWM}{C^{\infty}(M)}

\newcommand{\set}[1]{\left\{#1\right\}}

\newcommand{\frkd}{\mathfrak d}

\newcommand{\frkg}{\mathfrak g}
\newcommand{\frkh}{\mathfrak h}

\newcommand{\frkt}{\mathfrak t}

\newcommand{\frky}{\mathfrak y}

\newcommand{\frkX}{\mathfrak X}

\def\qed{\hfill ~\vrule height6pt width6pt depth0pt}

\newcommand{\pairing}[1]{\left\langle #1\right\rangle}

\newcommand{\pibracket}[1]{\left [ #1\right ]_{\pi}}
\newcommand{\conpairing}[1]{\left\langle  #1\right\rangle }

\newcommand{\jet}{\mathfrak{J}}

\newcommand{\jetd}{\mathbbm{d}}
\newcommand{\dev}{\mathfrak{D}}

\newcommand{\pie}{^\prime}
\newcommand{\Id}{\mathbf{1}}

\newcommand{\e}{\mathbbm{e}}
\newcommand{\p}{\mathbbm{p}}
\newcommand{\id}{\mathbbm{i}}
\newcommand{\jd}{\mathbbm{j}}

\newcommand{\dM}{\mathrm{d}}

\newcommand{\djet}{\mathrm{d}_{\jet}}
\newcommand{\dimension}{\mathrm{dim}}

\newcommand{\Hom}{\mathrm{Hom}}
\newcommand{\Der}{\mathrm{Der}}
\newcommand{\Inn}{\mathrm{Inn}}
\newcommand{\Out}{\mathrm{Out}}

\newcommand{\gl}{\mathfrak {gl}}

\newcommand{\Ker}{\mathrm{Ker}}

\newcommand{\ad}{\mathrm{ad}}

\newcommand{\Img}{\mathrm{Im}}

\begin{document}
\title{On Deformations of Lie Algebroids}

\author[Yunhe Sheng]{Yunhe Sheng}
\address{Mathematics School $\&$ Institute of Jilin University,\\
Changchun 130012, China\\
AND\\
 Department of Mathematics, Dalian University of
Technology,\\
 Dalian 116023,  China}

\email{shengyh@jlu.edu.cn}

\thanks{Research partially supported by NSF of China
(10871007), US-China CMR Noncommutative Geometry
(10911120391/A0109), China Postdoctoral Science Foundation
(20090451267), Science Research Foundation for Excellent Young
Teachers of  Mathematics School of Jilin University.}
\date{}



\begin{abstract}
For any Lie algebroid $A$, its $1$-jet bundle $\jet A$ is a Lie
algebroid naturally and there is a representation $\pi:\jet
A\longrightarrow\dev A$. Denote by $\dM_\jet$ the corresponding
coboundary operator. In this paper, we realize the deformation
cohomology of a Lie algebroid $A$ introduced by M. Crainic and I.
Moerdijk as the cohomology of a subcomplex
$(\Gamma(\Hom(\wedge^\bullet\jet A,A)_{\dev A}),\dM_\jet)$ of the
cochain complex $(\Gamma(\Hom(\wedge^\bullet\jet A,A)),\dM_\jet)$.
\end{abstract}
\maketitle

\section{Introduction}
The notion of Lie algebroid was introduced by Pradines in 1967, it
is a generalization of Lie algebras and tangent bundles. A Lie
algebroid over a manifold $M$ is a vector bundle $A\longrightarrow
M$ together with a Lie bracket $[\cdot,\cdot]$ on the section space
$\Gamma(A)$ and a bundle map $a:A\longrightarrow TM$, called the
anchor, satisfying the compatibility condition:
$$
[X,fY]=f[X,Y]+a(X)(f)Y,\quad\forall~X,~Y\in\Gamma(A),~f\in
C^\infty(M).
$$
We usually denote a Lie algebroid by $(A,[\cdot,\cdot],a)$, or $A$
if there is no confusion. See \cite{Fernandes,Mkz:GTGA} for more
details about Lie algebroids.

 In
\cite{Marius deformation cohomology}, Crainic and Moerdijk studied
the cohomology theory underlying deformations of Lie algebroids,
where they defined the deformation cohomology of a Lie algebroid
$(A,[\cdot,\cdot],a)$ and denote by
$\mathrm{H}^\bullet_{\mathrm{def}}(A)$. Any deformation of the Lie
bracket $[\cdot,\cdot]$  gives rise to a cohomology class in
$\mathrm{H}^2_{\mathrm{def}}(A)$. But in general, this cohomology
does not come from a representation of the Lie algebroid
$(A,[\cdot,\cdot],a)$. The deformation complex was also given by
Grabowska, Grabowski and Urba$\rm\acute{n}$ski in \cite{graboski},
where the authors studied Lie brackets on affine bundles.

 In general, there is  no natural
adjoint representation for Lie algebroids. For a Lie algebroid $A$,
the action of $\Gamma(A)$ on itself via the bracket is generally not
$C^\infty(M)$-linear in the first entry. There is a natural Lie
algebroid structure on the $1$-jet bundle $\jet A$. We call $\jet A$
the jet Lie algebroid of $A$. In fact, there is a natural Lie
algebroid structure on the $k$-jet bundle $\jet^k A$. There is also
a representation of the jet Lie algebroid $\jet A$ on $A$, which we
denote by $\pi:\jet A\longrightarrow\dev A$, where $\dev A$ is the
gauge Lie algebroid of $A$. This representation was first given in
\cite{secondcc}, where the authors call it the jet adjoint
representation of Lie algebroids. It was further studied by Blaom in
\cite{Blaum1,Blaum2}, where  the author call this representation the
adjoint representation of Lie algebroids. Recently, Camilo Arias
Abad and Marius Crainic suggest to define the adjoint representation
of Lie algebroids using representations up to homotopy \cite{Camilo
rep upto homotopy}. A similar notion under the name ``super
representations" was introduced independently by Alfonso Gracia-Saz
and Rajan Amit Mehta  in \cite{alfonso}.

The main purpose of this paper is to realize the deformation
cohomology of $A$ as some cohomology related to a representation. We
will see that the cohomology of the cochain complex
$(\Gamma(\Hom(\wedge^\bullet\jet A,A)_{\dev A}),\dM_\jet)$ is
isomorphic to the deformation cohomology, where $\dM_\jet$ is
decided by the representation $\pi:\jet A\longrightarrow\dev A$
(Theorem \ref{thm:main}).

The paper is organized as follows. In Section 2 we recall the
definition  of the deformation cohomology and we proved that for a
transitive Lie algebroid  $(A,[\cdot,\cdot],a)$, the deformation
cohomology is isomorphic to the cohomology of the Lie algebroid $A$
with coefficients in the adjoint representation. In Section 3 we
give the notion of the $k$-th differential operator bundle
$\Hom(\wedge^k\jet E,E)_{\dev E}$ and  establish the $k$-th
differential operator bundle sequence. In particular, for a Lie
algebroid $A$, we obtain a subcomplex
$(\Gamma(\Hom(\wedge^\bullet\jet A,A)_{\dev A}),\dM_\jet)$ of the
cochain complex $(\Gamma(\Hom(\wedge^\bullet\jet A,A)),\dM_\jet)$
associated with the representation $\pi:\jet A\longrightarrow\dev
A$. In Section 4 we prove that the cohomology of the subcomplex
$(\Gamma(\Hom(\wedge^\bullet\jet A,A)_{\dev A}),\dM_\jet)$ is
isomorphic to the deformation cohomology and we also give some
interesting examples.

\section{The deformation cohomology}
In this paper, $E\longrightarrow M$  is a vector bundle  with the
base manifold $M$. $\dM$ is the usual differential on forms. $\jetd$
is the coboundary operator associated with the complex
$(\Gamma(\Hom(\wedge^\bullet\dev E,E)_{\jet E}),\jetd)$. $\dM_\jet$
is the coboundary operator associated with the complex
$(\Gamma(\Hom(\wedge^\bullet\jet A,A)_{\dev A}),\dM_\jet)$.

 Recall that a multiderivation of degree $n$ of a
vector bundle $E$ is a skew-symmetric multi-linear map $
D:\Gamma(\wedge^{n}E)\longrightarrow\Gamma(E), $ such that for any
$f\in C^\infty(M)$ and $u_1,\cdots,u_n\in\Gamma(E)$, we have
 $$
D(u_1,\cdots,fu_n)=fD(u_1,\cdots,u_n)+\sigma_D(u_1,\cdots,u_{n-1})(f)u_n,
$$
where $ ~\sigma_D:\wedge^{n-1}E\longrightarrow TM $ is the symbol
of $D$. Denote by $\mathrm{D}^n(E)$ the set of multiderivations of
degree $n$. It is known that \cite{Marius deformation cohomology}
$\mathrm{D}^n(E)$ is the space of sections of a vector bundle
$\dev^nE\longrightarrow M$ which fits  into a short exact sequence
of vector bundles:
\begin{equation}\label{seq:der n}
0 \longrightarrow \Hom(\wedge^{n} E,
E)\stackrel{}{\longrightarrow}
                  {\dev^nE} \stackrel{}{\longrightarrow} \Hom(\wedge^{n-1}E,TM) \longrightarrow    0.
\end{equation}
In particular, $\dev^1E=\dev E$ is the gauge Lie algebroid of the
 frame bundle
 $\huaF(E)$, which is also called the covariant differential operator bundle of $E$ (see \cite[Example
3.3.4]{Mkz:GTGA}). The corresponding Atiyah sequence is as follows:
\begin{equation}\label{Seq:DE}
\xymatrix@C=0.5cm{0 \ar[r] & \gl(E)  \ar[rr]^{\id} &&
                \dev{E}  \ar[rr]^{\jd} && TM \ar[r]  & 0.
                }
\end{equation}

In \cite{Marius deformation cohomology}, the deformation complex of
a Lie algebroid $(A,[\cdot,\cdot],a)$ is defined as the complex
$(C^\bullet_{\mathrm{def}}(A),\delta)$ in which the $n$-cochains
$D\in C^\bullet_{\mathrm{def}}(A)$ are multi-linear skew-symmetric
maps
$$
D:\Gamma(\wedge^nA)\longrightarrow\Gamma(A),
$$
which are multiderivations and the coboundary operator is given by
\begin{eqnarray*}
\delta(D)(u_0,\cdots,u_n)&=&\sum_i(-1)^i[u_i,D(u_0,\cdots,\widehat{u_i},\cdots,u_n)]\\
&&+\sum_{i<j}(-1)^{i+j}D([u_i,u_j],u_0,\cdots,\widehat{u_i},\cdots,\widehat{u_j},\cdots,u_n).
\end{eqnarray*}
The deformation cohomology of a Lie algebroid $A$, denote by
$\mathrm{H}_{\mathrm{def}}^\bullet(A)$, is  the cohomology of the
cochain complex
$(C^\bullet_{\mathrm{def}}(A),\delta)$.\vspace{3mm}

Recall that a  Lie algebroid $(A,[\cdot,\cdot],a)$ is called a
transitive Lie algebroid if the anchor $a:A\longrightarrow TM$ is
surjective, and we have the following exact sequence of Lie
algebroids:
\begin{equation}
\xymatrix@C=0.5cm{0 \ar[r] & L \ar[rr]^{} &&
                A  \ar[rr]^{a} && TM \ar[r]  & 0,
                }
\end{equation}
where, $L=\ker(a)$ is a  bundle of Lie algebras. For transitive
Lie algebroids, there is a well defined adjoint representation
$\ad:A\longrightarrow\dev L$ of the Lie algebroid $A$ on the
vector bundle $L$, which is given by
$$
\ad_uX=[u,X],\quad\forall ~u\in\Gamma(A),~X\in\Gamma(L).
$$
Let $\dM_A$ be the coboundary operator  associated with the
adjoint representation. Denote the corresponding cochain complex
by $(C^\bullet(A;\ad),\dM_A)=(\Gamma(\Hom(\wedge^\bullet
A,L)),\dM_A)$, and the cohomology by $\mathrm{H}^\bullet(A;\ad)$.
In the following, we will show that $\mathrm{H}^\bullet(A;\ad)$ is
isomorphic to $\mathrm{H}^\bullet_{\mathrm{def}}(A)$. \vspace{3mm}

For any $n$-cochain $D\in C^n_{\mathrm{def}}(A)$, denote by
$D_a:\Gamma(\wedge^nA)\longrightarrow \frkX(M)$ the composition of
the anchor $a$ and the multiderivation $D$, i.e.
$$
D_a(u_1,\cdots,u_n)=a\big(
D(u_1,\cdots,u_n)\big),\quad\forall~u_1,\cdots,u_n\in\Gamma(A).
$$
Denote by $C^\bullet_a(A)$ the set of $D_a$, i.e.
\begin{equation}\label{eqn:ca}
C^n_a(A)=\{D_a\mid~\forall~D\in C^n_{\mathrm{def}}(A)\}.
\end{equation}
Obviously, for any $f\in C^\infty(M),$ we have
$$
D_a(u_1,\cdots,fu_n)=f
D_a(u_1,\cdots,u_n)+\sigma_D(u_1,\cdots,u_{n-1})(f)a(u_n),
$$
which implies that the symbol $\sigma_{D_a}$ of $D_a$ and the
symbol $\sigma_{D}$ are equal.

Furthermore, we can define the differential
$\delta:C^n_a(A)\longrightarrow C^{n+1}_a(A)$ by setting (we use
the same notation of the coboundary operator of the deformation
complex)
\begin{eqnarray*}
\nonumber\delta(D_a)(u_0,\cdots,u_n)&=&\sum_{i=0}^n(-1)^i[a(u_i),D_a(u_0,\cdots,\widehat{u_i},\cdots,u_n)]\\\label{d}
&&+\sum_{i<j}(-1)^{i+j}D_a([u_i,u_j],u_0,\cdots,\widehat{u_i},\cdots,\widehat{u_j},\cdots,u_n).
\end{eqnarray*}
Then we have
$$
\delta(D_a)=(\delta(D))_a.
$$
\begin{pro}\label{pro:complex CaA} With the above notations, for a
transitive Lie algebroid $A$, $C^n_a(A)$ which is defined by
\eqref{eqn:ca} is the space of sections of a vector bundle which we
denote by $\dev^n_aA$ and fits into the following exact sequence:
\begin{equation}\label{seq:Da}
0 \longrightarrow \Hom(\wedge^{n} A,
TM)\stackrel{}{\longrightarrow}
                  {\dev^n_aA} \stackrel{}{\longrightarrow} \Hom(\wedge^{n-1}A,TM) \longrightarrow    0.
\end{equation}
Furthermore, the complex $(C^\bullet_a(A),\delta)$ is acyclic.
\end{pro}
\pf The exact sequence (\ref{seq:Da}) follows from applying the
anchor $a$ to the exact sequence (\ref{seq:der n}). For any
multiderivation $D\in D^nA$, we have
$$
\sigma_{\delta(D)}=\delta(\sigma_D)+(-1)^{n+1}a\circ
D=\delta(\sigma_D)+(-1)^{n+1}D_a.
$$
If $\delta(D_a)=0$, we have $(\delta(D))_a=0$. Since
$\sigma_D=\sigma_{D_a}$, we have
$$
\sigma_{\delta(D)}=\sigma_{(\delta(D))_a}=0.
$$
Therefore, we obtain
$$D_a=(-1)^{n}\delta(\sigma_D)=(-1)^{n}\delta(\sigma_{D_a}),$$
which implies that $D_a$ is exact and this completes the proof.
\qed\vspace{3mm}

Therefore, we can obtain the following result which is also given
in \cite{Marius deformation cohomology}.

\begin{cor}
With the above notations, for a transitive Lie algebroid $A$, the
cohomology of $A$ with coefficients in the adjoint representation
is isomorphic to the deformation cohomology, i.e. we have
$$
\mathrm{H}^\bullet(A,\ad)\cong
\mathrm{H}^\bullet_{\mathrm{def}}(A).
$$
\end{cor}
\pf Obviously, the cochain complex $(C^\bullet(A;\ad),\dM_A)$ is a
subcomplex of $(C^\bullet_{\mathrm{def}}(A),\delta)$. For any $D\in
C^n_{\mathrm{def}}(A)$, write $D=D_a+D_L$, for some
$D_L\in\Gamma(\Hom(\wedge^nA,L))$. It is straightforward to see that
$$
\delta(D)=\delta(D_a)+\dM_A(D_L).
$$
Therefore, If $\delta(D)=0$, we have $\delta(D_a)=0$ and
$\dM_A(D_L)=0$. By Proposition \ref{pro:complex CaA}, the complex
$(C^\bullet_a(A),\delta)$ is acyclic, thus we have
$\mathrm{H}^\bullet(A,\ad)\cong
\mathrm{H}^\bullet_{\mathrm{def}}(A)$. \qed

\section{The complex $(\Hom(\wedge^k\jet A,A)_{\dev A},\dM_\jet)$}

The $1$-jet vector bundle ${\jet} E$ of the vector bundle $E$ (see
\cite{Jet bundle} for more details about jet bundles) is defined as
follows. For any $m\in M$, $({\jet}{E})_m$ is defined as a quotient
of local sections of $E$. Two local sections $u_1$ and $u_2$ are
equivalent and we denote this by $u_1\sim u_2$ if
$$ u_1(m)=u_2(m) ~~\mbox{ and }~~ \dM\pairing{u_{1
},\xi}_m=\dM\pairing{u_{2 },\xi}_m, \quad\forall~
\xi\in\Gamma(E^*).$$ So any $\mu\in ({\jet}{E})_m$ has a
representative $u\in\Gamma(E)$ such that
  $\mu=[u]_m$.
Let $\p$ be the projection which sends $[u]_m$ to $u(m)$.  Then
$\Ker\p \cong  \Hom(TM,E)$ and there is a short exact sequence,
called the jet sequence of $E$,
\begin{equation}\label{Seq:JetE}
\xymatrix@C=0.5cm{0 \ar[r] & \Hom(TM,E) \ar[rr]^{\quad\quad\quad\e}
&&
                {\jet}{E} \ar[rr]^{\p} && E \ar[r]  & 0.
                }
\end{equation}
from which it is straightforward to see that $\jet E$ is a finite
dimensional vector bundle. Moreover, $\Gamma(\jet{E})$ is isomorphic
to $\Gamma(E) \oplus \Gamma (T^*M \otimes E)$ as an $\Real$-vector
space.

In \cite{CLomni}, the authors proved that  the jet bundle $\jet E$
 may be
considered as an $E$-dual bundle of $\dev E$, i.e.
\begin{eqnarray*}  {\jet E}  &\cong &
\set{\nu\in \Hom( \dev{E}  ,E )\,|\,  \nu(\Phi)=\Phi\circ
\nu(\Id_E),\quad\forall ~~ \Phi\in \gl(E )}.
\end{eqnarray*}

A natural  nondegenerate symmetric $E$-valued pairing
$\conpairing{\cdot,\cdot}_E$ between $\jet E$ and $\dev{E}$ is
given by
\begin{eqnarray}\nonumber
\conpairing{\mu,\frkd}_E=\conpairing{\frkd,\mu}_E &\defbe& \frkd
u,\quad\forall ~~ \mu=[u]_m\in {\jet
E},~u\in\Gamma(E),~\frkd\in\dev{E}.
\end{eqnarray}
Moreover, this pairing  is $C^\infty(M) $-linear and satisfies the
following properties:
\begin{eqnarray*}
\conpairing{\mu, \Phi }_E &=& \Phi\circ \p(\mu),\quad\forall ~ ~\Phi\in \gl(E),~\mu\in{\jet E};\\
\conpairing{ {\frky} ,\frkd}_{E} &=& {\frky}\circ
\jd(\frkd),\quad\forall ~~ \frky\in \Hom(TM,E),~\frkd\in\dev{E}.
\end{eqnarray*}
For $k\geq 2$, the $k$-th skew-symmetric jet bundle
$\Hom(\wedge^k\dev E,E)_{\jet E}$ is defined in \cite{CLS2}:
\begin{eqnarray}\label{def:jet k}
\Hom(\wedge^k\dev E,E)_{\jet E}&\defbe& \set{\mu\in \Hom(\wedge^k
\dev E, E)~|~ \Img (\mu_\natural)\subset \jet E},
\end{eqnarray}
in which $\mu_\natural:\wedge^{k-1}\dev E\longrightarrow\Hom(\dev
E, E)$ is given by
\begin{equation}\label{eqn:natural map}
 \mu_\natural(\frkd_1,\cdots,
\frkd_{k-1})(\frkd_k)=\mu(\frkd_1,\cdots,
\frkd_{k-1},\frkd_k),\quad\forall ~\frkd_1,\cdots,\frkd_k\in\dev E.
\end{equation}
Furthermore, the authors proved that
$(\Gamma(\Hom(\wedge^\bullet\dev E,E)_{\jet E}),\jetd)$ is a
subcomplex of the cochain complex $(\Gamma(\Hom(\wedge^\bullet\dev
E,E),\jetd)$, where $\jetd$ is the coboundary operator associated
with the gauge Lie algebroid $\dev E$ with the obvious action on the
vector bundle $E$. In particular,
$\jetd:\Gamma(E)\longrightarrow\Gamma(\jet E)$ satisfies the
following formula which is very useful:
\begin{equation}\label{eqn:d fX}
\jetd(fX)=\dM f\otimes X+f\jetd X,\quad
\forall~X\in\Gamma(E),~f\in C^\infty(M).
\end{equation}
Furthermore, $\Gamma (\jet E)$ is an invariant subspace of the Lie
derivative $L_{\frkd}$ for any
 $\frkd \in\Gamma(\dev{E})$, which is defined  as follows:
\begin{eqnarray}\nonumber
\conpairing{L_{\frkd}\mu,\frkd\pie}_{E}&\defbe&
\frkd\conpairing{\mu,\frkd\pie}_{E}-\conpairing{\mu,[\frkd,\frkd\pie]_{\dev}}_{E},
\quad\forall~ \mu \in \Gamma(\jet{E}), ~
~\frkd\pie\in\Gamma(\dev{E}).
\end{eqnarray}

Considering the corresponding cohomology groups of the cochain
complex $(\Gamma(\Hom(\wedge^\bullet\dev E,E)_{\jet E}),\jetd)$, we
have
\begin{thm}\label{Thm:longExact}
For the cochain complex $C(E)=((\Gamma(\Hom(\wedge^\bullet\dev
E,E)_{\jet E}),\jetd)$, we have $\mathrm{H}^k(C(E))=0$, for all
$k=0,1,2,\cdots$. In other words, there is a long exact sequence:
\begin{eqnarray*}
  0\lon\Gamma(E)\stackrel{\jetd}{\longrightarrow}\Gamma(\jet{E})
\stackrel{\jetd}{\longrightarrow} \Gamma(\Hom(\wedge^2\dev
E,E)_{\jet
E})\stackrel{\jetd}{\longrightarrow}\cdots\\
\stackrel{\jetd}{\longrightarrow}
\Gamma(\Hom(\wedge^n\dev E,E)_{\jet E})\lon 0,
\end{eqnarray*}
where $n=\dimension M+1$.
\end{thm}
By this theorem,  the authors studied the deformation of omni-Lie
algebroids as well as its automorphism groups in
\cite{CLS2}.\vspace{3mm}

Assume that, for  the moment, the rank of $E$ is $r\geq2$. Similar
to (\ref{def:jet k}), we can define the $k$-th differential operator
bundle $\Hom(\wedge^k\jet E,E)_{\dev E}$ by
$$
\Hom(\wedge^k\jet E,E)_{\dev E} \triangleq \set{\frkd\in
\Hom(\wedge^k \jet E, E)~|~ \Img (\frkd_\natural)\subset \dev E},
\quad (k\geq 2).
$$
Next we study the property of the bundle $\Hom(\wedge^k\jet
E,E)_{\dev E}$ and give its corresponding exact sequence of vector
bundles.

\begin{pro}\label{Pro:muk=phi gamma}
For any $\frkd\in \Hom(\wedge^k\jet E,E)_{\dev E},$ there is a
unique bundle map $\lambda_{\frkd}\in \Hom(\wedge^{k-1}E,TM)$ such
that for any $\eta\in \Hom(TM,E),  ~\mu_i\in \jet{E}$, we have
\begin{equation}\label{eqn:j k dk} \frkd(
\mu_1\wedge\cdots\wedge\mu_{k-1}\wedge \eta)= \eta\circ
{\lambda_{\frkd}}(\p\mu_1\wedge\cdots\wedge\p\mu_{k-1}).
\end{equation}
\end{pro}
\pf For any $\alpha\otimes u,~\beta\otimes v\in\Hom(TM,E)$, where
$\alpha,~\beta\in\Omega^1(M)$ and $u,~v\in\Gamma(E)$, since $\frkd$
is skew-symmetric, we have
\begin{eqnarray*}
\frkd(\mu_1,\cdots,\mu_{k-2},\alpha\otimes u,\beta\otimes
v)&=&\langle\jd\circ\frkd_\natural(\mu_1,\cdots,\mu_{k-2},\alpha\otimes
u),\beta\rangle v\\
&=&-\langle\jd\circ\frkd_\natural(\mu_1,\cdots,\mu_{k-2},\beta\otimes
v),\alpha\rangle u,
\end{eqnarray*}
where the notation $\cdot_\natural$ is given by (\ref{eqn:natural
map}). Since the rank of $E$ was assumed to be bigger than one, it
follows that
$$
\jd\circ\frkd_\natural(\mu_1,\cdots,\mu_{k-2},\eta)=0,\quad\forall
~\eta\in\Hom(TM,E).
$$
Therefore, $\jd\circ\frkd_\natural$ factors though $\p$, i.e.
there is a unique $\lambda_{\frkd}\in\Hom(\wedge^{k-1}E,TM)$ such
that
$$
\jd\circ \frkd_\natural(\mu_1\wedge\cdots\wedge\mu_{k-1})=
{\lambda_{\frkd}}(\p\mu_1\wedge\cdots\wedge\p\mu_{k-1}),
$$
which yields the conclusion.  \qed \vspace{3mm}

We will write $\jd(\frkd)=\lambda_{\frkd}$ by (\ref{eqn:j k dk}).
For any $\Phi\in\Hom(\wedge^kE,E)$, $~\id(\Phi)\in\Hom(\wedge^k\jet
E,E)_{\dev E} $ is given by
\begin{equation}\label{eqn:i k}
\id(\Phi)(\mu_1,\cdots,\mu_{k})=\Phi(\p\mu_1,\cdots,\p\mu_{k}).
\end{equation}

\begin{thm}
For any $k\geq1$, we have the following exact sequence:
\begin{equation}\label{seq:1}
0 \longrightarrow \Hom(\wedge^k E,
E)\stackrel{{\id}}{\longrightarrow}
                  {\Hom(\wedge^k\jet E,E)_{\dev E}} \stackrel{\jd}{\longrightarrow} \Hom(\wedge^{k-1}E,TM) \longrightarrow    0.
\end{equation}
\end{thm}
\pf For any $\lambda\in \Hom(\wedge^{k-1}E,TM) $, define
$\widehat{\lambda}\in\Hom(\wedge^k\jet E,E)_{\dev E}$  by
$$
\widehat{\lambda}(\mu_1,\cdots,\mu_k)=\sum_{i}(-1)^{i+1}(\mu_i-\gamma\p\mu_i)(\lambda(\p\mu_1,\cdots,\widehat{\p\mu_i},
\cdots,\p\mu_k)),
$$
for any split $\gamma:E\longrightarrow\jet E$ of (\ref{Seq:JetE}).
Evidently, $\mu_i-\gamma\p\mu_i\in\Hom(TM,E)$ and
$\widehat{\lambda}\in\Hom(\wedge^k\jet E,E)$. Furthermore, for any
$\eta\in\Hom(TM,E)$, we have
$$
\widehat{\lambda}(\mu_1,\cdots,\mu_{k-1},\eta)=(-1)^{k+1}\eta\circ\lambda(\p\mu_1,\cdots,
\p\mu_{k-1}),
$$
which means that $\widehat{\lambda}\in\Hom(\wedge^k\jet E,E)_{\dev
E}$ and $\jd((-1)^{k+1}\widehat{\lambda})=\lambda$, i.e.
 the bundle map $\jd$ is surjective. By (\ref{eqn:i
k}), the definition of $\id$, it is obvious that $\jd\circ \id=0$.

If $\frkd\in\Hom(\wedge^k\jet E,E)_{\dev E}$ satisfies
$\jd(\frkd)=0$, we have
$$
\frkd(\mu_1,\cdots,\mu_{k-1},\eta)=\eta\circ\jd(\frkd)(\p\mu_1,\cdots,
\p\mu_{k-1})=0,
$$
which implies $\frkd$ factors through $\p$, i.e. there is a unique
$\Phi\in\Hom(\wedge^kE,E)$ such that
$$
\frkd(\mu_1,\cdots,\mu_{k})=\Phi(\p\mu_1,\cdots,\p\mu_{k}).
$$
This completes the proof of the exactness of (\ref{seq:1}). \qed
\vspace{3mm}

We call exact sequence (\ref{seq:1}) the $k$-th differential
operator bundle sequence.
\begin{rmk}
If the rank of the vector bundle is $r=1$, when $n\geq2$, we should
extend the definition of $\Hom(\wedge^n\jet E,E)_{\dev E}$ to
satisfy the exact sequence \eqref{seq:1}.
\end{rmk}

Associated with any Lie algebroid $(A,[\cdot,\cdot],a)$, there is
a bundle map $\pi:\jet A\longrightarrow\dev A$ which is given by
\cite{CLomni}
\begin{equation}\label{eqn:def pi}
\pi(\jetd u)(v)=[u,v],\quad\forall~u,~v\in\Gamma(A),
\end{equation}
and a bracket $[\cdot,\cdot]_\pi$ on $\Gamma(\jet A)$ by setting
\begin{equation}\label{eqn:bracket pi}
\pibracket{\mu,\nu}\defbe L_{\pi(\mu)}\nu-L_{\pi(\nu)}\mu-
\jetd\conpairing{\pi(\mu),\nu}_A=L_{\pi(\mu)}\nu-i_{\pi(\nu)}\jetd\mu.
\end{equation}
It turns out that $(\jet A,[\cdot,\cdot]_\pi,\jd\circ\pi)$ is a Lie
algebroid together with the representation $\pi$. We give a list of
several useful formulas here which will be used later. The proof is
straightforward by (\ref{eqn:d fX}), (\ref{eqn:def pi}),
(\ref{eqn:bracket pi}) and we leave it to  readers.
\begin{lem}\label{formulas}
For any $u,~v\in \Gamma(A),~\omega,~\theta\in\Omega^1(M),~f\in
\CWM$, we have
\begin{eqnarray} \label{eqn:1}[\jetd u,\jetd
v]_\pi&=&\jetd[u,v],\\~[\jetd u,\omega\otimes v]_\pi&=&L_{a(u)}\omega\otimes v+\omega\otimes [u,v],\\
~[\omega\otimes u,\theta\otimes
v]_\pi&=&\conpairing{a(u),\theta}\omega\otimes
v-\conpairing{a(v),\omega}\theta\otimes u,\\\label{eqn:2}
L_\frkd(\dM f\otimes v)&=&\dM f\otimes \frkd
v+\dM(\jd(\frkd)f)\otimes v,\\\label{eqn:3} \pi(\dM f\otimes
v)(u)&=&-a(u)(f)v,
\end{eqnarray}
and
\begin{eqnarray}\label{eqn:4}
 \jd(\pi(\jetd u))=a(u).
\end{eqnarray}
\end{lem}

For more information about (\ref{eqn:bracket pi}), see
\cite{CLomni}. Denote by $\djet$ the associated coboundary operator
in the cochain complex $(\Gamma(\Hom(\wedge^\bullet\jet
A,A)),\djet)$. Furthermore, in \cite{CLS2}, by using the theory of
Manin pairs, the authors proved that $ (\dev A,\jet A)$ is an
$A$-Lie bialgebroid. Therefore, $(\Gamma(\Hom(\wedge^\bullet\jet
A,A)_{\dev A}),\djet)$ is a subcomplex of the cochain complex
$(\Gamma(\Hom(\wedge^\bullet\jet A,A)),\djet)$. In fact, we have

\begin{pro}
For any $\frkd\in\Hom(\wedge^k\jet A,A)_{\dev A}$, we have
$$\jd(\djet\frkd)=\delta(\jd(\frkd))+(-1)^{k+1}a\circ\frkd\circ\jetd.$$
More precisely, for any $u_1,\cdots,u_k\in\Gamma(A)$, we have
\begin{equation}\label{eqn:jd d d}
\jd(\djet\frkd)(u_1,\cdots,u_k)=\delta(\jd(\frkd))(u_1,\cdots,u_k)+(-1)^{k+1}a\circ\frkd(\jetd
u_1,\cdots,\jetd u_k),
\end{equation}
where $\delta $ is given by \eqref{d}.
\end{pro}
\pf 
For any $\mu_1,\cdots,\mu_k,\dM f\otimes v\in\Gamma(\jet A)$, we
have
\begin{eqnarray*}
&&\djet\frkd(\mu_1,\cdots,\mu_k,\dM f\otimes
v)\\&=&\sum_{i=1}^k(-1)^{i+1}\pi(\mu_i)\frkd(\mu_1,\cdots,\widehat{\mu_i},\cdots,\mu_k,\dM
f\otimes v)\\
&&+(-1)^k\pi(\dM f\otimes v)\frkd(\mu_1,\cdots,\mu_k)\\
&&+\sum_{i<j\leq
k}(-1)^{i+j}\frkd([\mu_i,\mu_j]_\pi,\mu_1,\cdots,\widehat{\mu_i},\cdots,\widehat{\mu_j},\cdots,\mu_k,\dM
f\otimes v)\\
&&+\sum_{i }(-1)^{i+k+1}\frkd([\mu_i,\dM f\otimes
v]_\pi,\mu_1,\cdots,\widehat{\mu_i},\cdots,\mu_k).
\end{eqnarray*}
By straightforward computations, we have
\begin{eqnarray*}
&&\sum_{i=1}^k(-1)^{i+1}\pi(\mu_i)\frkd(\mu_1,\cdots,\widehat{\mu_i},\cdots,\mu_k,\dM
f\otimes
v)\\
&=&\sum_{i=1}^k(-1)^{i+1}\pi(\mu_i)\big(\jd(\frkd)(\p\mu_1,\cdots,\widehat{\p\mu_i},\cdots,\p\mu_k)(f)v\big)\\
&=&\sum_{i=1}^k(-1)^{i+1}\jd(\frkd)(\p\mu_1,\cdots,\widehat{\p\mu_i},\cdots,\p\mu_k)(f)\pi(\mu_i)(v)\\
&&+\sum_{i=1}^k(-1)^{i+1}\jd(\pi(\mu_i))(\jd(\frkd)(\p\mu_1,\cdots,\widehat{\p\mu_i},\cdots,\p\mu_k)(f))v.
\end{eqnarray*}
By \eqref{eqn:3}, we obtain
$$
(-1)^k\pi(\dM f\otimes
v)\frkd(\mu_1,\cdots,\mu_k)=(-1)^{k+1}a\circ\frkd(\mu_1,\cdots,\mu_k)(f)v.
$$
It is obvious that
\begin{eqnarray*}
&&\sum_{i<j\leq
k}(-1)^{i+j}\frkd([\mu_i,\mu_j]_\pi,\mu_1,\cdots,\widehat{\mu_i},\cdots,\widehat{\mu_j},\cdots,\mu_k,\dM
f\otimes v)\\
&=&\sum_{i<j\leq
k}(-1)^{i+j}\jd(\frkd)(\p[\mu_i,\mu_j]_\pi,\p\mu_1,\cdots,\widehat{\p\mu_i},\cdots,\widehat{\p\mu_j},\cdots,\p\mu_k)(f)v.
\end{eqnarray*}
By \eqref{eqn:2}, we have
\begin{eqnarray*}
[\mu_i,\dM f\otimes v]_\pi&=&L_{\pi(\mu_i)}(\dM f\otimes
v)-i_{\pi(\dM
f\otimes v)}\jetd\mu_i\\
&=&\dM f\otimes \pi(\mu_i)(v)+\dM\circ\jd(\pi(\mu_i))(f)\otimes
v-i_{\pi(\dM f\otimes v)}\jetd\mu_i.
\end{eqnarray*}
Consequently, we have
\begin{eqnarray*}
&&\sum_{i }(-1)^{i+k+1}\frkd([\mu_i,\dM f\otimes
v]_\pi,\mu_1,\cdots,\widehat{\mu_i},\cdots,\mu_k)\\
&=&\sum_{i
}(-1)^{i+k+1}\Big((-1)^{k-1}\frkd(\mu_1,\cdots,\widehat{\mu_i},\cdots,\mu_k,\dM f\otimes\pi(\mu_i)(v))\\
&&+(-1)^{k-1}\frkd(\mu_1,\cdots,\widehat{\mu_i},\cdots,\mu_k,\dM\circ\jd(\pi(\mu_i))(f)\otimes
v)\\
&&+(-1)^{k}\frkd(\mu_1,\cdots,\widehat{\mu_i},\cdots,\mu_k,i_{\pi(\dM
f\otimes v)}\jetd\mu_i)\Big)\\
&=&\sum_i(-1)^i\jd(\frkd)(\p\mu_1,\cdots,\widehat{\p\mu_i},\cdots,\p\mu_k)
(f)\pi(\mu_i)(v)\\
&&+\sum_i(-1)^i\jd(\frkd)(\p\mu_1,\cdots,\widehat{\p\mu_i},\cdots,\p\mu_k)(\jd(\pi(\mu_i))(f))v\\
&&+\sum_i(-1)^{i+1}\conpairing{\frkd(\mu_1,\cdots,\widehat{\mu_i},\cdots,\mu_k),i_{\pi(\dM
f\otimes v)}\jetd\mu_i}_A.
\end{eqnarray*}
Therefore, we have
\begin{eqnarray*}
&&\djet\frkd(\mu_1,\cdots,\mu_k,\dM f\otimes
v)\\
&=&(-1)^{k+1}a\circ\frkd(\mu_1,\cdots,\mu_k)(f)v\\
&&+\sum_{i<j\leq
k}(-1)^{i+j}\jd(\frkd)(\p[\mu_i,\mu_j],\p\mu_1,\cdots,\widehat{\p\mu_i},\cdots,\widehat{\p\mu_j},
\cdots,\p\mu_k)(f)v\\
&&+\sum_i(-1)^{i+1}[\jd(\pi(\mu_i)),\jd(\frkd)(\p\mu_1,\cdots,\widehat{\p\mu_i},\cdots,\p\mu_k)](f)v\\
&&+\sum_i(-1)^{i+1}\conpairing{\frkd(\mu_1,\cdots,\widehat{\mu_i},\cdots,\mu_k),i_{\pi(\dM
f\otimes v)}\jetd\mu_i}_A.
\end{eqnarray*}
Now assume that $\mu_1=\jetd u_1,\cdots, \mu_k=\jetd u_k$, we get
\begin{eqnarray*}
&&\djet\frkd(\jetd u_1,\cdots, \jetd u_k,\dM f\otimes
v)\\
&=&(-1)^{k+1}a\circ\frkd(\jetd u_1,\cdots, \jetd u_k)(f)v\\
&&+\sum_{i<j\leq
k}(-1)^{i+j}\jd(\frkd)([u_i,u_j],u_1,\cdots,\widehat{u_i},\cdots,\widehat{u_j},
\cdots,u_k)(f)v\\
&&+\sum_i(-1)^{i+1}[a(u_i),\jd(\frkd)(u_1,\cdots,\widehat{u_i},\cdots,u_k)](f)v\\
&=&\Big(\delta(\jd(\frkd))(u_1,\cdots,u_k)+(-1)^{k+1}a\circ\frkd(\jetd
u_1,\cdots,\jetd u_k)\Big)(f)v.
\end{eqnarray*}
This implies that
$$
\jd(\djet\frkd)=\delta(\jd(\frkd))+(-1)^{k+1}a\circ\frkd\circ\jetd.
\qed
$$

\section{Infinitesimal Deformations of Lie algebroids}

Denote by $\mathrm{H}^\bullet(\jet A;A)$ the resulting cohomology of
$(\Gamma(\Hom(\wedge^\bullet\jet A,A)_{\dev A}),\djet)$. The main
result in this section is
\begin{thm}\label{thm:main}
Given a Lie algebroid $A$, the cohomology $\mathrm{H}^\bullet(\jet
A;A)$ is isomorphic to the deformation cohomology of $A$, i.e.
$$
\mathrm{H}^\bullet(\jet A;A)\cong
\mathrm{H}_{\mathrm{def}}^\bullet(A).
$$
\end{thm}
\pf First we prove that there is a one-to-one correspondence between
$\Gamma(\Hom(\wedge^k\jet A,A)_{\dev A})$ and
$C^k_{\mathrm{def}}(A)$ at the level of cochains. For any
$\frkd\in\Gamma(\Hom(\wedge^k\jet A,A)_{\dev A})$, define
$D_{\frkd}\in C^k_{\mathrm{def}}(A)$ by
\begin{equation}\label{eqn:def D mu k}
D_{\frkd}(u_1,\cdots,u_k)\triangleq\frkd(\jetd u_1,\cdots,\jetd
u_k).
\end{equation}
Follow from
\begin{eqnarray*}
D_{\frkd}(u_1,\cdots,fu_k)&=&\frkd(\jetd u_1,\cdots,\jetd
(fu_k))\\
&=&f\frkd(\jetd u_1,\cdots,\jetd u_k)+\frkd(\jetd u_1,\cdots,\dM
f\otimes u_k)\\
&=&fD_{\frkd}(u_1,\cdots,u_k)+\jd(\frkd)(u_1,\cdots,u_{k-1})(f)u_k,
\end{eqnarray*}
we know that $D_{\frkd}$ is well defined and the following equality
holds:
$$
\sigma_{D_{\frkd}}=\jd(\frkd).
$$
Conversely, for any $D\in C^k_{\mathrm{def}}(A)$, define
$\frkd_D\in\Gamma(\Hom(\wedge^k\jet A,A)_{\dev A})$ by
\begin{eqnarray*}
\frkd_D(\jetd u_1,\cdots,\jetd
u_k)&=&D(u_1,\cdots,u_k),\end{eqnarray*} and
\begin{eqnarray*}
\frkd_D(\jetd
u_1,\cdots,\dM f\otimes u_k)&=&\sigma_D(u_1,\cdots,u_{k-1})(f)u_k.
\end{eqnarray*}
By (\ref{eqn:d fX}), it is straightforward to see that $\frkd_D$
is well defined and satisfies
$$
\jd(\frkd_D)=\sigma_D.
$$
Furthermore, obviously we have
$$
\frkd_{D_{\frkd}}=\frkd,\quad D_{\frkd_D}=D,
$$
which implies that, at the level of cochains, there is a one-to-one
correspondence between $\Gamma(\Hom(\wedge^k\jet A,A)_{\dev A})$ and
$C^k_{\mathrm{def}}(A)$.

If $\frkd$ is closed, i.e. $\djet\frkd=0$, then follows from
$\pi(\jetd u)(v)=[u,v]$ and $\jetd[u,v]=[\jetd u,\jetd v]_\pi$, we
have
\begin{eqnarray*}
&&\delta(D_{\frkd})(u_0,\cdots,u_k)\\
&=&\sum_i(-1)^i[u_i,D_{\frkd}(u_0,\cdots,\widehat{u_i},\cdots,u_k)]\\
&&+\sum_{i<j}(-1)^{i+j}D_{\frkd}([u_i,u_j],u_0,\cdots,\widehat{u_i},\cdots,\widehat{u_j},\cdots,u_k))\\
&=&\sum_i(-1)^i[u_i,\frkd(\jetd u_0,\cdots,\widehat{\jetd u_i},\cdots,\jetd u_k)]\\
&&+\sum_{i<j}(-1)^{i+j}\frkd(\jetd[u_i,u_j],\jetd
u_0,\cdots,\widehat{\jetd u_i},
\cdots,\widehat{\jetd u_j},\cdots,\jetd u_k))\\
&=&\djet\frkd(\jetd u_0,\cdots,\jetd u_k) \\
&=&0.
\end{eqnarray*}
If  $\frkd$ is exact, i.e. there is some $\frkt\in
\Hom(\wedge^{k-1}\jet A,A)_{\dev A}$ such that $\frkd=\djet\frkt$,
then we have
\begin{eqnarray*}
&&D_{\frkd}(u_1,\cdots,u_k)\\&=&\frkd(\jetd u_1,\cdots,\jetd u_k)\\
&=&(\djet\frkt)(\jetd u_1,\cdots,\jetd u_k)\\
&=&\sum_i(-1)^{i+1}[u_i,\frkt(\jetd u_1,\cdots,\widehat{\jetd u_i},\cdots,\jetd u_k)]\\
&&+\sum_{i<j}(-1)^{i+j}\frkt(\jetd[u_i,u_j],\jetd
u_1,\cdots,\widehat{\jetd u_i},
\cdots,\widehat{\jetd u_j},\cdots,\jetd u_k)\\
&=&\delta(D_{\frkt})(u_1,\cdots,u_k).
\end{eqnarray*}
Conversely, if $D\in C^k_{\mathrm{def}}(A)$ is closed, first we
have
\begin{eqnarray*}
(\djet\frkd_D)(\jetd u_0,\cdots,\jetd
u_k)=\delta(D)(u_0,\cdots,u_k)=0.
\end{eqnarray*}
Then for any $m\leq k$, any $f_l,~l=m,\cdots,k$, by (\ref{eqn:d
fX}), we have
\begin{eqnarray*}
&&(\djet\frkd_D)(\jetd u_0,\cdots,\jetd u_{m-1},\dM f_m\otimes
u_m,\cdots,\dM f_k\otimes u_k)\\&=&(\djet\frkd_D)(\jetd
u_0,\cdots,\jetd u_{m-1},\jetd (f_mu_m),\cdots,\jetd
(f_ku_k))\\&&-f_m\cdots f_k(\djet\frkd_D)(\jetd u_0,\cdots,\jetd
u_{m-1},\jetd u_m,\cdots,\jetd u_k)\\
&=&0,
\end{eqnarray*}
which implies that $\frkd_D$ is closed. Similarly, if $D$ is
exact, $\frkd_D$ is also exact. The proof of the theorem is
completed. \qed \vspace{3mm}

\begin{cor}
 $\frkd\in\Gamma(\Hom(\wedge^k\jet A,A)_{\dev A})$ is closed if and only of
 $\frkd\mid_{\jetd\Gamma(A)}$ is closed, i.e. for any
 $u_0,\cdots,u_k$,
$$
(\djet\frkd)(\jetd u_0,\cdots,\jetd u_k)=0.
$$
\end{cor}

Next we examine the cohomology $\mathrm{H}^\bullet(\jet A;A)$ in
low degrees. \vspace{3mm}

$\bullet$ In degree 0,  $u\in\Gamma(A)$ is  closed means that $u$
belongs to the center $Z(\Gamma(A))$ of the infinite-dimensional Lie
algebra $\Gamma(A)$, i.e.
$$
\mathrm{H}^0(\jet A;A)=Z(\Gamma(A)).
$$
In fact,  for any $\mu\in\jet A$, $\djet u(\mu)=0$ is equivalent to
the condition that for any $v\in\Gamma(A),~\omega\in\Omega^1(M)$,
$$\djet u(\jetd v)=0,\quad\djet u(\omega\otimes v)=0.$$
On the other hand, we have
\begin{eqnarray*}
\djet u(\jetd v)&=&\pi(\jetd v)(u)=[v,u],
\end{eqnarray*}
and
\begin{eqnarray*}
\djet u(\omega\otimes v)&=&\pi(\omega\otimes
v)(u)=-\langle\omega,a(u)\rangle v.
\end{eqnarray*}
Thus we have
$$
\djet u=0\Longleftrightarrow \left\{\begin{array}{c}[u,v]=0,~\forall~v\in\Gamma(A)\\
a(u)=0.\end{array}\right.
$$
However, if $[u,v]=0,$ for any $~v\in\Gamma(A)$, then for any $f\in
\CWM$, we have $[u,fv]=0$, which implies that $a(u)(f)=0,~$ for any
$f\in \CWM$. This happens exactly when $a(u)=0$. Thus we have
$$
\djet u=0\Longleftrightarrow u\in Z(\Gamma(A)).
$$

$\bullet$ In degree 1,  $\frkd\in\Gamma(\dev A)$ is closed if and
only if $\frkd\in\Der(A)$, where $\Der(A)$ denotes the set of
derivatives of the Lie algebroid $A$. In fact,  $\frkd\in\Gamma(\dev
A)$ is closed if and only if  for any
$u,~v\in\Gamma(A),~\omega,~\theta\in\Omega^1(M)$, the following
equalities hold:
$$
(\djet\frkd)(\jetd u,\jetd v)=0,\quad (\djet\frkd)(\jetd
u,\theta\otimes v)=0,\quad(\djet\frkd)(\omega\otimes u,\theta\otimes
v)=0.
$$

On the other hand, by Lemma \ref{formulas}, we have
\begin{eqnarray*}
&&(\djet\frkd)(\jetd u,\jetd v)\\
&=&\pi(\jetd
u)\conpairing{\frkd,\jetd v}_A-\pi(\jetd v)\conpairing{\frkd,\jetd
u}_A-\conpairing{\frkd,[\jetd
u,\jetd v]_\pi}_A\\
&=&[u,\frkd(v)]-[v,\frkd(u)]-\frkd([u,v]),\\
&&(\djet\frkd)(\jetd u,\theta\otimes v)\\&=&\pi(\jetd
u)\conpairing{\frkd,\theta\otimes v}_A-\pi(\theta\otimes
v)\conpairing{\frkd,\jetd u}_A-\conpairing{\frkd,[\jetd
u,\theta\otimes v]_\pi}_A\\
&=&[u,\conpairing{\theta,\jd\frkd}v]+\conpairing{\theta,a(\frkd(u))}v-\conpairing{\jd\frkd,L_{a(u)}\theta}v-\conpairing{\theta,\jd\frkd}[u,v]\\
&=&\big(\conpairing{\theta,a(\frkd(u))}+\conpairing{\theta,[a(u),\jd\frkd]}\big)v,
\end{eqnarray*}
and
\begin{eqnarray*}
&&(\djet\frkd)(\omega\otimes u,\theta\otimes
v)\\&=&\pi(\omega\otimes u)\conpairing{\frkd,\theta\otimes
v}_A-\pi(\theta\otimes v)\conpairing{\frkd,\omega\otimes
u}_A-\conpairing{\frkd,[\omega\otimes u,\theta\otimes v]_\pi}_A\\
&=&-\conpairing{\jd\frkd,\theta}\conpairing{a(v),\omega}u+\conpairing{\jd\frkd,\omega}\conpairing{a(u),\theta}v\\&&-\conpairing{\frkd,\conpairing{a(u),\theta}\omega\otimes
v-\conpairing{a(v),\omega}\theta\otimes u}_A\\
&=&0.
\end{eqnarray*}
Therefore, $(\djet\frkd)(\jetd u,\jetd v)=0$ if and only if $\frkd$
is a derivation with respect to the Lie bracket, i.e.
$\frkd\in\Der(A)$. Furthermore, $\frkd\in\Der(A)$ also implies that
$(\djet\frkd)(\jetd u,\theta\otimes v)=0$. This follows from the
next lemma.

\begin{lem}
If $\frkd\in\Der(A)$, i.e. $\frkd$ is a derivation with respect to
the Lie bracket of $A$, then we have
$$
a(\frkd(u))+[a(u),\jd\frkd]=0,\quad\forall~u\in\Gamma(A).
$$
\end{lem}
\pf Since $\frkd$ is a derivation, for any $f\in \CWM$, we have
$$
\frkd[u,fv]=[\frkd(u),fv]+[u,\frkd(fv)].
$$
Furthermore, we have
\begin{eqnarray*}
\frkd[u,fv]&=&\frkd(f[u,v]+a(u)(f)v)\\
&=&\jd\frkd(f)[u,v]+f\frkd[u,v]+a(u)(f)\frkd(v)+\jd\frkd(a(u)(f))v,
\end{eqnarray*}
and
\begin{eqnarray*}
~[\frkd(u),fv]+[u,\frkd(fv)]&=&f[\frkd(u),v]+a(\frkd(u))(f)v+a(u)(f)\frkd(v)\\
&&+f[u,\frkd (v)]+\jd\frkd(f)[u,v]+a(u)\jd\frkd(f)v.
\end{eqnarray*}
Thus we have
 $$
a(\frkd(u))+[a(u),\jd\frkd]=0,\quad\forall~u\in\Gamma(A).~\qed
$$

Therefore, $\djet\frkd=0$ if and only if $\frkd\in\Der(A)$. If
$\frkd$ is exact, i.e. $\frkd=\djet u$ for some $u\in\Gamma(A)$, we
have
\begin{eqnarray*}
\frkd(\jetd v)&=&\conpairing{\djet u,\jetd v}_A=\pi(\jetd
v)(u)=-[u,v],
\end{eqnarray*}
and
\begin{eqnarray*}
\frkd(\omega\otimes v)&=&\conpairing{\djet u,\omega\otimes
v}_A=\pi(\omega\otimes v)(u)=-\conpairing{a(u),\omega}v,
\end{eqnarray*}
which implies that $\frkd=-\ad_u$. Thus we have
$$
\mathrm{H}^1(\jet A;A)=\Out(A)=\Der(A)/\Inn(A).
$$

Let $(A,[\cdot,\cdot],a)$ be a fix Lie algebroid over the base
manifold $M$ and $I\subset \mathbb R$ be an integral. A 1-parameter
infinitesimal deformation of the Lie algebroid $A$ over $I$ is a
collection $A_t$ of Lie algebroids $A_t=(A,[\cdot,\cdot]_t,a_t)$
varying smoothly with respect to $t$. In \cite{Marius deformation
cohomology}, for a deformation $A_t=(A,[\cdot,\cdot]_t,a_t)$ of the
Lie algebroid $A$, the authors proved that there is an associated
2-cocycle $c_0\in C^2_{\mathrm{def}}(A)$ which is defined by
$$
c_0(X,Y)=\frac{\dM}{\dM t}[X,Y]_t\Big|_{t=0},\quad \forall~X,Y\in
\Gamma(A).
$$
In fact only a 2-cocycle can not contain all the information of
the deformation, this 2-cocycle should also define a Lie
bracket.\vspace{3mm}

Next we
 consider the 1-parameter infinitesimal deformation of the Lie algebroid $A$ of the following form:
\begin{equation}\label{eqn:deformation}
[X,Y]_t=[X,Y]+t\frkd(\jetd X,\jetd Y), \quad \forall~X,Y\in
\Gamma(A),
\end{equation}
where $\frkd\in\Gamma(\Hom(\wedge^2\jet A,A)_{\dev A})$. It is easy
to see that the anchors vary as follows:
$$
a_t=a+t\jd(\frkd).
$$
Since $[ \cdot, \cdot]_t$ should satisfy the Jacobi identity, we
can obtain
\begin{eqnarray}
\label{eqn:111}\frkd(\jetd[X,Y],\jetd Z)+[\frkd(\jetd X,\jetd
Y),Z]+c.p.&=&0,
\end{eqnarray}
and
\begin{eqnarray}
\label{eqn:222} D_{\frkd}(D_{\frkd}(X,Y),Z)+c.p.&=&0,
\end{eqnarray}
which implies that $\djet\frkd=0$ and $D_{\frkd}$ (see (\ref{eqn:def
D mu k})) itself defines a Lie bracket.\vspace{3mm}

We summarize the discussion in the following proposition.
\begin{pro}
For any $1$-parameter infinitesimal deformation of the Lie algebroid
$A$ of the form \eqref{eqn:deformation},
$\frkd\in\Gamma(\Hom(\wedge^2\jet A,A)_{\dev A})$ is a $2$-cocycle
such that $D_{\frkd}:\Gamma(\wedge^2A)\longrightarrow\Gamma(A)$
defines a Lie algebroid structure on  $A$.
\end{pro}

\begin{rmk}
In general, we can study higher order deformations and versa
deformations, see \cite{Fialowski1,Fialowski2} for more details.
However, for our objects we leave this study for later
consideration.
\end{rmk}

By (\ref{eqn:def pi}), we can write (\ref{eqn:deformation}) as
\begin{eqnarray*}
[X,Y]_t&=&\pi(\jetd X)(Y)+t\frkd(\jetd X)(Y)\\
&=&(\pi+t\frkd)(\jetd X)(Y).
\end{eqnarray*}
Thus, for any $t\in I$, $[\cdot,\cdot]_t$ is a Lie bracket iff
$$[\pi+t\frkd,\pi+t\frkd]=0,$$ which holds if and only if
$$[\pi,\frkd]=0,\quad [\frkd,\frkd]=0.$$ It is straightforward to
see that it is equivalent to (\ref{eqn:111}) and (\ref{eqn:222}). In
particular, if we only condition the deformation in the following
form:
\begin{equation}\label{1}
[X,Y]_\frkd=[X,Y]+\frkd(\jetd X,\jetd Y),
\end{equation}
obviously we have
\begin{thm}
With the above notations, \eqref{1} defines a deformation of the Lie
algebroid $(A,[\cdot,\cdot],a)$ for some
$\frkd\in\Gamma(\Hom(\wedge^2\jet A,A)_{\dev A})$  if and only if
 $\frkd$ satisfies the Maurer-Cartan equation:
\begin{equation}\label{maurer-cartan}
\dM_\jet\frkd+\frac{1}{2}[\frkd\wedge\frkd]=0.
\end{equation}
\end{thm}
\begin{rmk}
In fact, for any $\frkd\in\Gamma(\Hom(\wedge^2\jet A,A)_{\dev A})$,
if we consider the graph $\huaG_\frkd$ which is given by
$$
\huaG_\frkd=\{\frkd(\mu)+\mu|~\forall~\mu\in\jet A\}\subset\dev
A\oplus \jet A,
$$
\eqref{maurer-cartan} also means that  $\huaG_\frkd$ is a Dirac
structure. See \cite{CLS1} and  Theorem 7.8 in  \cite{CLS2} for more
details about Dirac structures.
\end{rmk}
\begin{ex}\label{ex:g}{\em
If the Lie algebroid $A$ is a Lie algebra $\frkg$, we have $\jet
\frkg=\frkg,~\dev \frkg=\gl(\frkg)$ and $\Hom(\wedge^2 \jet
\frkg,\frkg)_{\dev \frkg}=\Hom(\frkg\wedge\frkg,\frkg) $. The
resulting cohomology $\mathrm{H}^\bullet(\jet \frkg,\frkg)$ turns
out to be the cohomology of the Lie algebra $\frkg$ with
coefficients in the adjoint representation. See
\cite{Fialowski1,Fialowski2} for more details.}
\end{ex}
\begin{ex}{\em
If the Lie algebroid $A$ is the tangent Lie algebroid $TM$, we have
already known that all the deformations are trivial \cite{Marius
deformation cohomology}. In fact, in this case, it is evident that
the gauge Lie algebroid $\dev (TM)$ is isomorphic to the jet Lie
algebroid $\jet (TM)$. Therefore, the cohomology of cochain complex
$(\Gamma(\Hom(\wedge^\bullet\jet ({TM}),TM)_{\dev(TM)}),\djet)$ is
isomorphic to the cohomology of cochain complex
$(\Gamma(\Hom(\wedge^\bullet\dev ({TM}),TM)_{\jet(TM)}),\jetd)$. By
Theorem \ref{Thm:longExact}, we know that all the deformations are
trivial. It also implies that the tangent Lie algebroid is rigid.}
\end{ex}

\begin{ex}\rm{
We consider the deformation of Lie algebroid $(A,[\cdot,\cdot],a)$
by a two cocycle $\frkd=\djet N$, where
$N\in\Gamma(\gl(A))\subset\Gamma(\dev A)$. For all
$u,v\in\Gamma(A)$, we have
\begin{eqnarray*}
\frkd(\jetd u,\jetd v)=\djet N(\jetd u,\jetd v)&=&[u,Nv]+[Nu,v]-N[u,v]\\
&=&[u,v]_N.
\end{eqnarray*}
If $N$ is a Nijenhuis operator, we can obtain the deformation of
$A$ as follows:
\begin{equation} [u,v]_t=[u,v]+t[u,v]_N,\quad
a_t=a+ta\circ N,\quad\forall~u,v\in\Gamma(A).
\end{equation}
The deformation of the Lie algebroid $A$ by a Nijenhuis operator
is trivial, i.e. $\Id+tN$ is an isomorphism from  the Lie
algebroid $(A,[\cdot,\cdot]_t,a_t)$ to the Lie algebroid
$(A,[\cdot,\cdot],a)$.

In particular, if $A$ is the cotangent bundle Lie algebroid of some
Poisson manifold, we can consider the compatibility condition of
Poisson structures and Nijenhuis structures, i.e. Poisson-Nijenhuis
structures. For more information about Poisson-Nijenhuis structures
on oriented 3D-manifolds, see \cite{chensheng}.}
\end{ex}
\begin{ex}\rm{
For a Lie algebra $\frkg$, there is a Lie-Poisson structure
$\pi_1$ on $\frkg^*$, we can consider the deformation of the
corresponding Lie algebroid by a quadratic Poisson structure
$\pi_2$. The corresponding 2-cocycle $\Omega_{\pi_2}$ is given by
$$
\Omega_{\pi_2}(\xi,\eta)=L_{\pi_2^\sharp(\xi)}\eta-L_{\pi_2^\sharp(\eta)}\xi-\dM\pi_2(\xi,\eta),\quad\forall~\xi,~\eta\in\Omega^1(\mathbb
\frkg^*).
$$
Obviously, $\Omega_{\pi_2}$ defines a Lie bracket if and only if
$\pi_2$ is a Poisson structure, $\Omega_{\pi_2}$ is closed is
equivalent to the condition that $\pi_2$ and $\pi_1$ are
compatible:
$$
[\pi_1,\pi_2]=0.
$$
For more information about quadratic deformation of Lie-Poisson
structures on $\mathbb R^3$, see \cite{LLS}}.
\end{ex}
\begin{ex}\rm{
We consider a special deformation of a 4-dimensional Lie algebra
$\frkh$, which is the direct sum of a 3-dimensional Lie subalgebra
$\frkg$ and a 1-dimensional center $\mathbb Re$. As shown in
Example \ref{ex:g}, $\mathrm{H}^\bullet(\jet\frkh,\frkh)$ is just
the Lie algebra cohomology of $\frkh$ with coefficient in the
adjoint representation. For any $D\in\Der(\frkg)$,
$\Omega_D:\wedge^2\frkh\longrightarrow\frkh$ is given by£º
\begin{equation}\label{eqn:omegaD}
\Omega_D(X,Y)=0;\quad\Omega_D(X,e)=De,\quad\forall~X,Y\in\frkg.
\end{equation}
Obviously, $\Omega_D$ defines a Lie bracket and since $e$ is a
center of $\frkh$, we have $\Omega_D$ is closed if and only if $D$
is a derivation of $\frkg$. Therefore, this problem can also be
considered as the extension of a 3-dimensional Lie algebra by a
derivation. In \cite{sheng}, the author gives a classification of
such extensions using Poisson geometry method and therefore obtains
the classification of 4-dimensional Lie algebras at the end (see
also \cite{Fialowski3,Fialowski4}).}
\end{ex}
\subsection*{Acknowledgment}We would like to give our special thanks to
 Zhuo Chen and Zhangju Liu for
their helpful suggestions and comments. We also would like to give
our warmest thanks to Chenchang Zhu for the help during we stayed in
Courant Research Center, G$\ddot{\rm{o}}$ttingen, where a part of
work was done. Special thanks are  given to referees for very useful
comments.


\begin{thebibliography}{999}
\bibitem{Camilo rep upto homotopy}
C. Arias Abad and M. Crainic, Representations up to homotopy of
Lie algebroids, arXiv:0901.0319.

\bibitem{Blaum1}
A. Blaom, Geometric structures as deformed infinitesimal symmetries,
\emph{Trans. Amer. Math. Soc.}, 358 (2006), 3651-3671.

\bibitem{Blaum2}
 A.  Blaom, Lie algebroids and Cartan's method of equivalence,
arXiv:math/0509071v3.


\bibitem{chensheng}
B. Chen and Y. Sheng, Poisson-Nijenhuis structures on oriented
3D-manifolds.  {\em Rep. Math. Phys.} 61 (3) (2008), 361-380.

\bibitem{CLomni}
Z. Chen and Z. Liu, Omni-Lie algebroids, \emph{J. Geom. Phys.} 60
(2010), 799-808.

\bibitem{CLS1}
Z. Chen, Z. Liu and Y. Sheng, Dirac structures of omni-Lie
algebroids, arXiv:0802.3819., to appear in \emph{International J.
Math.}

\bibitem{CLS2}
Z. Chen, Z. Liu and Y. Sheng, $E$-Courant algebroids, \emph{Int.
Math. Res. Not.} Vol. 2010, No. 22, pp. 4334-4376.
\bibitem{secondcc}
M. Crainic and R. L. Fernandes,  Secondary characteristic classes of
Lie algebroids. \emph{Quantum field theory and noncommutative
geometry,} 157--176, Lecture Notes in Phys., 662, Springer, Berlin,
2005.

\bibitem{Marius deformation cohomology}
M. Crainic and I. Moerdijk, Deformation of Lie brackets:
cohomological aspects, \emph{J. Eur. Math. Soc.} (JEMS) 10 (2008),
no. 4, 1037-1059.


\bibitem{Fernandes}
 R. Fernandes, Lie algebroids, holonomy and characteristic classes.
\emph{Adv. Math.} 170 (2002), no. 1, 119--179.

\bibitem{Fialowski1}
A. Fialowski, Deformations of Lie algebras, \emph{Math. USSR
Sbornik}, Vol. 55 (1986), No. 2, 467-473.

\bibitem{Fialowski2}
A. Fialowski, An example of formal deformations of Lie algebras,
\emph{NATO Conf. Proceedings,} Kluwer (1988), 375-401.

\bibitem{Fialowski3}
 A. Fialowski and M. Penkava, Deformations of Four Dimensional Lie
 Algebras, \emph{Comm. Contemp. Math.} 9 (2007), 41-79.

\bibitem{Fialowski4}
A. Fialowski and M. Penkava, Moduli spaces of low dimensional real
Lie algebras, \emph{J. Math. Phys.} 49 (2008), 073507.






\bibitem{graboski}
K. Grabowska, J. Grabowski and P. Urba$\rm\acute{¡än}$ski, Lie
brackets on affine bundles, \emph{Ann. Global Anal. Geom.} 24
(2003), 101-130.








\bibitem{alfonso}
A. Gracia-Saz and R.A. Mehta, Lie algebroid structures on double
vector bundles and representation theory of Lie algebroids, \emph{
Adv. Math.} 223 (2010), 1236-1275.


\bibitem{LLS}
Q. Lin, Z. Liu and Y. Sheng, Quadratic deformations of Lie-Poisson
structures on $\mathbb R^3$,  {\it Lett. Math. Phys.} 83 (2008),
217-229.





\bibitem{Mkz:GTGA}
K. Mackenzie, \emph{General theories of Lie groupoids and Lie
algebroids}, Cambridge University Press, 2005.




\bibitem{Jet bundle}
D. J. Saunders, \emph{The Geometry of Jet Bundles}, Cambridge
University Press, Cambridge, 1989.

\bibitem{sheng}
Y. Sheng, Linear Poisson structures on $\mathbb R^4$, {\em J. Geom.
Phys.} 57 (2007), 2398-2410.





\end{thebibliography}
\end{document}